\documentclass[oneside,british]{amsart}
\usepackage[T1]{fontenc}
\usepackage[utf8]{inputenc}
\usepackage{geometry}
\usepackage{babel}
\usepackage{amsbsy}
\usepackage{amstext}
\usepackage{amsthm}
\usepackage{amssymb}
\usepackage{graphicx}
\usepackage{booktabs}
\usepackage{subcaption}
\usepackage[unicode=true,pdfusetitle,
 bookmarks=true,bookmarksnumbered=false,bookmarksopen=false,
 breaklinks=false,pdfborder={0 0 1},backref=false,colorlinks=false]
 {hyperref}

\numberwithin{equation}{section}
\numberwithin{figure}{section}

\theoremstyle{plain}
\newtheorem{thm}{Theorem}
\numberwithin{thm}{section}
\theoremstyle{plain}

\theoremstyle{plain}

\theoremstyle{definition}

\theoremstyle{remark}
\newtheorem{rem}[thm]{Remark}

\DeclareMathOperator{\HH}{\mathrm{H}}
\DeclareMathOperator{\tr}{\mathrm{tr}}
\DeclareMathOperator{\hyp}{\mathrm{hyp}}
\DeclareMathOperator{\ext}{\mathrm{ext}}
\newcommand{\NN}{\mathrm{NN}}

\title[A user's guide to PINNs in geometric analysis:\\ lessons from the asymptotic Plateau problem]{A user's guide to PINNs in geometric analysis:\\ lessons from the asymptotic Plateau problem}

\author{Tancredi Schettini Gherardini}
\thanks{University of Bonn and Max Planck Institute for Mathematics, Bonn, Germany.
E-mail: \href{mailto:tsg@math.uni-bonn.de}{tsg@math.uni-bonn.de}.}

\begin{document}

\begin{abstract}
This proceedings contribution elaborates on the findings of arXiv:2605.26234v2:
a joint work with Marco Usula, where we introduced a machine learning
framework based on physics-informed neural networks (PINNs), aimed at
constructing near-minimal discs in hyperbolic space asymptotic to a
prescribed knot at infinity. We used this method to provide numerical evidence
for a conjecture of Joel Fine relating minimal surfaces in $\HH^{4}$ to the
coefficients of the HOMFLY polynomial. This is a
methodological companion to that paper, based on a presentation given at 
the 2026 edition of the workshop ``DANGER: Data, Numbers, and Geometry''.
Rather than reviewing the results, which are presented extensively in the preprint 
above, we discuss the two aspects of the framework which, in our experience, 
determined whether the method worked at all. First, the geometry of the problem
must be encoded in the architecture of the model, so that the boundary
condition and asymptotics at infinity hold exactly
for every value of the learnable parameters -- leaving us with a single-component
loss function; second, the evaluation of the PDE residual must be engineered 
with care to ensure that complete trainings can be performed in a reasonable time. 
On the latter point, we describe two implementation techniques which are not spelled
out in detail in the original paper: replacing nested reverse-mode automatic 
differentiation  with the forward propagation of second-order jets, and compiling 
the computational graph of the residual once instead of rebuilding it at
every optimisation step. Together, on identical hardware, these two changes reduce the cost of a training step by a
factor of roughly forty to fifty. We hope these methodological discussions
can be useful for researchers in differential geometry and geometric 
analysis who wish to deploy PINNs on problems of their own.
\end{abstract}

\maketitle

\section{Introduction}

The \emph{asymptotic Plateau problem} asks whether a closed submanifold
$M^{k}\subset S^{n}$ of the sphere at infinity of hyperbolic space
$\HH^{n+1}$ bounds a complete, properly immersed minimal submanifold of
$\HH^{n+1}$. Anderson's classical theorem \cite{AndersonPlateau}
guarantees the existence of an area-minimising solution, but says
nothing about the structure of the full space of solutions. For
\emph{surfaces} in $\HH^{4}$ bounding \emph{knots} $K\subset S^{3}$,
Fine \cite{FineKnots} proved that (for generic $K$) the minimal
surfaces of fixed genus $g$ and self-intersection number $d$ bounding
$K$ form a compact oriented $0$-dimensional moduli space when $g=0$,
whose signed count is a knot invariant, and explained why these counts
should be closely related to the coefficients of the HOMFLY polynomial
\cite{HOMFLY} of $K$; this precise conjectural correspondence was 
numerically tested and corroborated in \cite{SchettiniUsula2026}.

The paper presented a numerical framework,
based on \emph{physics-informed neural networks} (PINNs)
\cite{raissi2017physicsinformeddeeplearning}, which produces
near-minimal discs in $\HH^{4}$ bounding a prescribed knot,
locates their self-intersections, and computes the sign of each double
point. Testing the resulting self-intersection numbers against the
HOMFLY coefficients -- for the unknot, the torus knots $T(3,2)$,
$T(5,2)$, $T(4,3)$ and $T(5,3)$, the figure-eight, three-twist,
Stevedore and square knots, and various mirror images -- we found
agreement with Fine's conjecture in every case. We
refer to \cite{SchettiniUsula2026} for the geometric background, the
precise statements, and the results; the trained models and the full
implementation are publicly available\footnote{\url{https://github.com/Tancredi-Schettini-Gherardini/deep_plateau}}.

The present note, prepared for the proceedings of the ``DANGER: Data, Numbers, and Geometry'' 
2026 workshop, where this work was presented, deliberately takes a different angle. It
is written for geometers and analysts who are considering PINNs as a
tool, and its subject is not \emph{what} we computed but \emph{how}.
Neural networks are being used across geometric analysis and adjacent
fields
\cite{cortes2026machinelearningapproachnirenberg,platt2026nonuniquenesssymmetriesnirenbergproblem,Hirst2025AInstein,HalversonRuehle2023MetricFlows,hirst2026minimisingwillmoreenergyneural,ZhouYe2023MinimalSurfacePINN,Hashimoto2025PINNMinimal},
including as a component of computer-assisted proofs
\cite{gomez2019computer,wang2025discoveryunstablesingularities}, and in
our experience the difference between a PINN that converges to
$10^{-6}$ residuals in under an hour on a laptop and one that stalls,
or takes days, is not primarily the network or the optimiser: it lies
in choices that are usually left undocumented. Two of these choices
form the theses of this note.

\begin{enumerate}
\item \textbf{The geometry of the problem should be encoded in the
model, not learned.} Our model satisfies the boundary condition at
infinity \emph{exactly, for every value of the learnable parameters},
and is \emph{asymptotically minimal by construction}: the analytically
known leading behaviour of solutions near the boundary is built into
the architecture, and the network only learns the interior. This
removes the boundary term from the loss entirely, and with it the
notoriously delicate balancing of boundary against interior penalties
(\S\ref{sec:geometry}); it yields a single-component loss encoding only
one geometric property, which in our case is \textit{minimality}.
\item \textbf{The evaluation of the PDE residual should be engineered,
not delegated.} A PINN loss for a second-order equation needs first and
second derivatives of the model with respect to its \emph{inputs} at
every collocation point, inside an outer optimisation loop that
differentiates with respect to the \emph{parameters}. The default
implementation -- nested reverse-mode automatic differentiation,
rebuilt at every step -- is wasteful in a way that compounds
multiplicatively. Choosing the direction of differentiation by counting
dimensions (forward for inputs, reverse for parameters), propagating
second-order jets through the network in closed form, and compiling the
resulting computational graph \emph{once}, sped up our training by a
factor of $40$--$50$ on identical hardware (\S\ref{sec:efficiency});
neither technique appears in the companion paper.
\end{enumerate}

Section \ref{sec:setting} recalls the minimal amount of geometric setup
needed to make this note self-contained and explains where the standard
PINN recipe struggles on this problem. Sections \ref{sec:geometry} and
\ref{sec:efficiency} develop the two theses. Section
\ref{sec:guidelines} distils the discussion into a short list of
guidelines which, we believe, transfer to many geometric variational
problems beyond the asymptotic Plateau problem.

\subsubsection*{Acknowledgements}
The results surveyed here were obtained jointly with Marco Usula, whom
the author thanks for the collaboration and for many conversations
reflected in this note. We also thank Ed Hirst for the invitation to the
DANGER workshop, and all the organisers for putting together such a great
event. Finally, we also thank once more Joel Fine for his
feedback and encouragement during the preparation of the main preprint.
We acknowledge the contribution of general-purpose AI agents (Claude Opus 4.8 and Fable 5)
in the implementation of the code and in the preparation of the first draft of this paper.
We acknowledge the support of the 2024
Max Planck-Humboldt Research Award, bestowed on Geordie Williamson by
the Max Planck Society and the Alexander von Humboldt Foundation and
hosted by Catharina Stroppel at the University of Bonn.

\section{The problem as a PINN task}\label{sec:setting}

\subsection{Minimal discs in half-space coordinates}\label{subsec:halfspace}

We work in the upper half-space model of hyperbolic space: on
$\mathbb{R}_{>0}\times\mathbb{R}^{n}$ with coordinates
$(X,\boldsymbol{Y})=(X,Y_{1},\dots,Y_{n})$, the hyperbolic metric is
$g_{\hyp}=X^{-2}\left(dX^{2}+|d\boldsymbol{Y}|^{2}\right)$, and the
boundary at infinity of the closure is the locus
$\{X=0\}\cong\mathbb{R}^{n}$, a stereographic chart of $S^{n}$ minus a
point. The object we seek is a
map
\[
u=(X,\boldsymbol{Y}):D^{2}\longrightarrow\overline{\HH}^{n+1}
\]
from the closed unit disc, mapping the interior to the interior, whose
boundary restriction $u_{|\partial D^{2}}$ parametrises a prescribed
embedded curve $\gamma:S^{1}\to\mathbb{R}^{n}$ (for $n=3$: a knot) in
the boundary at infinity, and whose interior image is minimal. In the
language of \cite{SchettiniUsula2026,UsulaBiharmonic, UsulaIsometricEmbeddings}, $u$ should be a
minimal \emph{p-immersion}. Minimality is expressed by the vanishing of
the tension field $\tau(u)=\tr_{u^{*}g_{\hyp}}\nabla du$, computed with
respect to the pull-back metric; in half-space coordinates, writing
$g=X^{-2}J^{\top}J$ for the pull-back metric ($J$ the Jacobian of $u$
in the disc coordinates), the components of $\tau(u)$ in the orthonormal
frame $(X\partial_{X},X\partial_{Y_{k}})$ read
\begin{align}
\tau^{X}(u) & =\frac{1}{X}\Bigl[\Delta_{g}X+\frac{1}{X}\Bigl(\sum_{k=1}^{n}|dY_{k}|_{g}^{2}-|dX|_{g}^{2}\Bigr)\Bigr],\label{eq:tau_X}\\
\tau^{Y_{k}}(u) & =\frac{1}{X}\Bigl[\Delta_{g}Y_{k}-\frac{2}{X}\,\langle dX,dY_{k}\rangle_{g}\Bigr],\qquad k=1,\dots,n,\label{eq:tau_Y}
\end{align}
where $\Delta_{g}$ is the Laplace--Beltrami operator of $g$, acting on
a scalar $f$ by
\begin{equation}
\Delta_{g}f=g^{ab}\partial_{ab}^{2}f+(\partial_{a}g^{ab})\partial_{b}f+g^{ab}(\partial_{a}\log\sqrt{\det g})\partial_{b}f.\label{eq:LB}
\end{equation}
The equation $\tau(u)=0$ is a quasi-linear second-order elliptic
system, degenerate at the boundary, with no general method available
to solve it explicitly. Training minimises the Monte Carlo $L^{2}$
norm of the residual,
\begin{equation}
\mathcal{L}(\theta)=\frac{1}{N}\sum_{i=1}^{N}\left|\tau(u_{\theta})\right|^{2}(p_{i}),\label{eq:loss}
\end{equation}
over a sample $\{p_{i}\}\subset D^{2}$ of interior \emph{collocation
points}; here $\theta$ denotes the learnable parameters of the model
$u_{\theta}$, and the pointwise norm is taken with respect to the
hyperbolic metric, i.e.
$|\tau|^{2}=(\tau^{X})^{2}+\sum_{k}(\tau^{Y_{k}})^{2}$ in the frame
above.

\subsection{The standard recipe, and where it strains}\label{subsec:vanilla}

The textbook PINN prescription for a boundary value problem
$\mathcal{F}[u]=0$ in $\Omega$, $\mathcal{B}[u]=0$ on $\partial\Omega$,
is to take $u_{\theta}$ to be a multi-layer perceptron (MLP) and to
minimise a weighted sum of an interior and a boundary penalty,
\[
\mathcal{L}(\theta)=\frac{1}{N_{\Omega}}\sum_{i}\bigl|\mathcal{F}[u_{\theta}](x_{i})\bigr|^{2}+\frac{\lambda}{N_{\partial\Omega}}\sum_{j}\bigl|\mathcal{B}[u_{\theta}](x_{j}')\bigr|^{2},
\]
the derivatives in $\mathcal{F}$ being supplied by automatic
differentiation. This recipe is general, and its very generality is
its weakness on a problem like ours, for three reasons. First, the
\emph{weighting problem}: $\lambda$ trades boundary accuracy against
interior accuracy, the optimal trade-off drifts during training, and a
wrong choice produces either a surface that misses its boundary curve
or one that fits the boundary while violating the equation; a sizeable
literature on adaptive weighting schemes exists precisely because this
balance is fragile. Second, the \emph{boundary is at infinity}: the
ambient metric blows up as $X\to0$, so uniform penalties on the map do
not translate into uniform geometric control, and a network that is
free to move the boundary even slightly produces maps that are not
asymptotic to the prescribed knot at all. Third, the \emph{regularity
theory is known and is not exploited}: by a theorem of Marx-Kuo
\cite{MarxKuoRenormalized}, a minimal p-submanifold of
$\overline{\HH}^{n+1}$ with smooth boundary at infinity meets the
boundary \emph{orthogonally}, and admits a precise polyhomogeneous
expansion there. A generic soft-constrained network knows nothing of
this and must spend its capacity discovering it.

The framework of \cite{SchettiniUsula2026} addresses these three points, as we now explain.

\section{Encoding the geometry in the model}\label{sec:geometry}

\subsection{A hard-constrained ansatz}\label{subsec:ansatz}

Instead of asking an MLP to represent $u$ directly, the model wraps the
network inside a composite map whose structure carries the geometry.
Fix an embedding $\gamma:S^{1}\to\mathbb{R}^{3}$ (the knot), a
\emph{boundary defining function} $\rho$ for the disc (a smooth
function on $D^2$, positive in the interior, vanishing on
$\partial D^{2}$, but whose exterior derivative is non-zero on $\partial D^{2}$), an \emph{extension operator} $\ext$ producing a map
$\ext(\gamma):D^{2}\to\mathbb{R}^{3}$ with
$\ext(\gamma)_{|\partial D^{2}}=\gamma$, and a \emph{decay exponent}
$k\in\{1,2\}$. The model is
\begin{equation}
u_{\theta}=\Bigl(\;\rho\,e^{\NN^{X}},\;\ext(\gamma)+\rho^{k}\,\NN^{\boldsymbol{Y}}\Bigr),\label{eq:ansatz}
\end{equation}
where $\NN=(\NN^{X},\NN^{\boldsymbol{Y}}):\mathbb{R}^{2}\to\mathbb{R}^{4}$
is the only learnable ingredient --- in all our experiments a plain
MLP with four hidden layers of width $64$ and $\tanh$ activations,
$12\,932$ parameters in total. Three structural facts hold for
\emph{every} value of $\theta$:
\begin{enumerate}
\item the first component is positive in the interior and vanishes
exactly on $\partial D^{2}$, so the image lies in $\HH^{4}$ and reaches
the boundary at infinity precisely over $\partial D^{2}$;
\item the boundary restriction of $u_{\theta}$ is exactly $\gamma$: the
factor $\rho^{k}$ kills the network correction at the boundary, where
$\ext(\gamma)$ restricts to $\gamma$;
\item consequently the loss \eqref{eq:loss} contains \emph{no boundary
term at all}: the weighting problem of \S\ref{subsec:vanilla}
disappears, because the constraint manifold has been parametrised away
rather than penalised.
\end{enumerate}
A fourth fact concerns initialisation, and we single it out because it
guided every later refinement: \emph{at $\theta=0$ the model reduces to
the explicit map $(\rho,\ext(\gamma))$}. The network is therefore best
understood as a \emph{learnable perturbation} of a built-in starting
point, and choosing $\rho$ and $\ext$ amounts to choosing the best
explicit approximation to a solution that one can write down. For the
round unknot $\gamma(\varphi)=(\cos\varphi,\sin\varphi,0)$ (we write
$(r,\varphi)$ for polar coordinates on $D^{2}$, reserving $\theta$ for
the learnable parameters) the best choice is obvious: the totally
geodesic copy of $\HH^{2}$,
\begin{equation}
u_{\HH^{2}}(r,\varphi)=\Bigl(\frac{1-r^{2}}{1+r^{2}},\,\frac{2r\cos\varphi}{1+r^{2}},\,\frac{2r\sin\varphi}{1+r^{2}},\,0\Bigr),\label{eq:geodesic}
\end{equation}
is an exact solution. Reading off its first component motivates the
\emph{stereographic} boundary defining function
$\rho_{\mathrm{st}}=(1-r^{2})/(1+r^{2})$, which we use throughout: with
it (and the extension below), $\theta=0$ reproduces
\eqref{eq:geodesic} \emph{exactly} when $\gamma$ is the round unknot.

\subsection{Choosing the extension: a regularity and asymptotics
story}\label{subsec:extension}

The extension operator is the heart of the construction, and finding
the right one took several iterations over the life of the project. We
record the sequence, because each step encodes a transferable lesson.

\subsubsection*{Radial fattening}
The earliest implementations extended the knot as a cone,
$\ext(\gamma)(r,\varphi)=c(r)\,\gamma(\varphi)$ with $c(r)=r$ or
$\frac{2r}{1+r^{2}}$, and smoothed variants. The latter restricts to
$\gamma$ at $r=1$ and, paired with $\rho_{\mathrm{st}}$, is again exact
for the round unknot; but for a generic curve it is merely continuous
at the origin, where the cone has a genuine singularity. Since the
network correction in \eqref{eq:ansatz} is smooth, no value of $\theta$
can repair a non-smooth baked-in ingredient: in experiments, training
stalls with the error concentrated near the origin, and ad hoc
smoothings of the cone tip remove the singularity only at the cost of
the boundary asymptotics. The simple lesson:
the regularity of every hard-coded ingredient is a ceiling on
the regularity of the model -- hard constraints transfer their
defects to the solution just as reliably as their virtues.

\subsubsection*{Harmonic extension}
Every smooth $\gamma:S^{1}\to\mathbb{R}^{3}$ has a unique harmonic
extension $\Gamma$ to the disc, which is smooth up to the boundary and
computable to machine precision from the Fourier coefficients of
$\gamma$: writing $z=x+iy$ and
$\gamma\sim A_{0}+\sum_{m\ge1}(A_{m}\cos m\varphi+B_{m}\sin m\varphi)$,
one has
$\Gamma=A_{0}+\sum_{m\ge1}\left(U_{m}A_{m}+V_{m}B_{m}\right)$ with
$U_{m}+iV_{m}=z^{m}$. The \emph{stereoharmonic} extension
$\ext(\gamma)=2\Gamma/(1+r^{2})$ again reproduces
\eqref{eq:geodesic} for the round unknot, and is now smooth for every
$\gamma$. This already trains well.

\subsubsection*{Biharmonic extension and orthogonality at
infinity}
The remaining defect concerns asymptotics rather than smoothness. As
recalled in \S\ref{subsec:vanilla}, an actual minimal p-submanifold
meets the boundary at infinity orthogonally; equivalently, it is
\emph{asymptotically minimal}, i.e.\ $|\tau|^{2}\to0$ at the boundary.
For a map of the form \eqref{eq:ansatz} one computes that orthogonality
holds if and only if the normal derivative
$\partial_{r}\bigl(\ext(\gamma)+\rho^{k}\NN^{\boldsymbol{Y}}\bigr)$ is
tangent to the curve along $\partial D^{2}$
\cite[\S3.2]{SchettiniUsula2026}. If $k=1$ the network can, in
principle, learn to correct a non-orthogonal extension; if $k=2$ its
correction decays too fast to do so, and the extension alone must carry
the right asymptotics. This suggests imposing a \emph{Neumann}
condition on the extension in addition to the Dirichlet one --- which
is exactly what the bi-Laplacian allows. The \emph{stereobiharmonic}
extension is $\ext(\gamma)=2\Gamma/(1+r^{2})$, where $\Gamma$ is the
unique solution of
\[
\Delta^{2}\Gamma=0,\qquad\Gamma_{|S^{1}}=\gamma,\qquad\partial_{r}\Gamma_{|S^{1}}=\gamma;
\]
mode by mode,
\[
\Gamma=A_{0}\,\tfrac{1+r^{2}}{2}+\sum_{m\ge1}W_{m}(r^{2})\left(U_{m}A_{m}+V_{m}B_{m}\right),\qquad W_{m}(r^{2})=\tfrac{m+1}{2}+\tfrac{1-m}{2}\,r^{2},
\]
so it is exactly as computable as the harmonic one. Paired with
$k=2$, the resulting model is orthogonal to the boundary --- hence
asymptotically minimal --- \emph{for every value of the learnable
parameters}, while still reducing to \eqref{eq:geodesic} at $\theta=0$
for the round unknot. This is the production configuration used for
all results in \cite{SchettiniUsula2026}.

\begin{figure}[t]
\centering
\begin{subfigure}[b]{0.44\textwidth}
\centering
\includegraphics[width=1\textwidth]{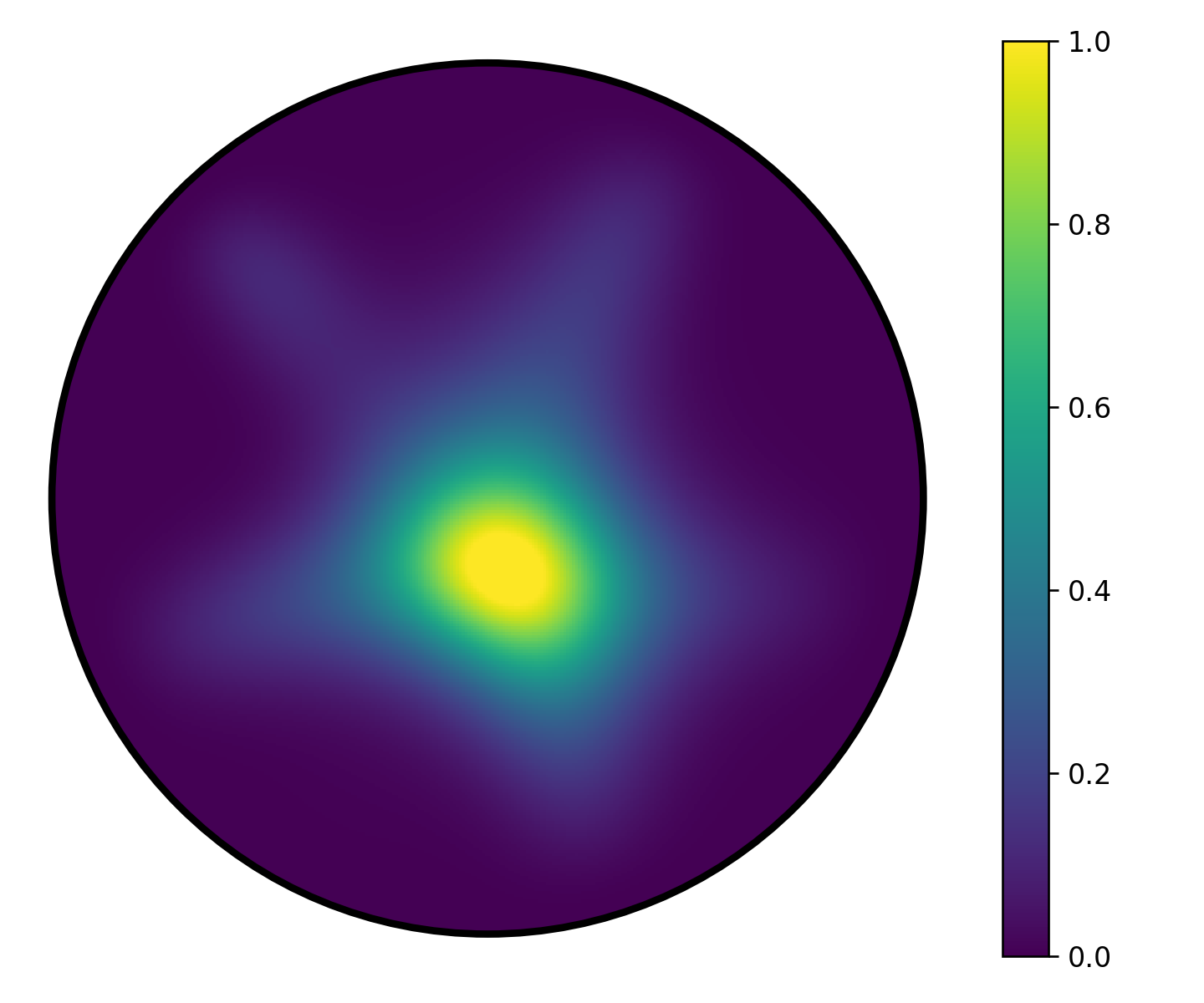}
\end{subfigure}
\hfill{}
\begin{subfigure}[b]{0.44\textwidth}
\centering
\includegraphics[width=1\textwidth]{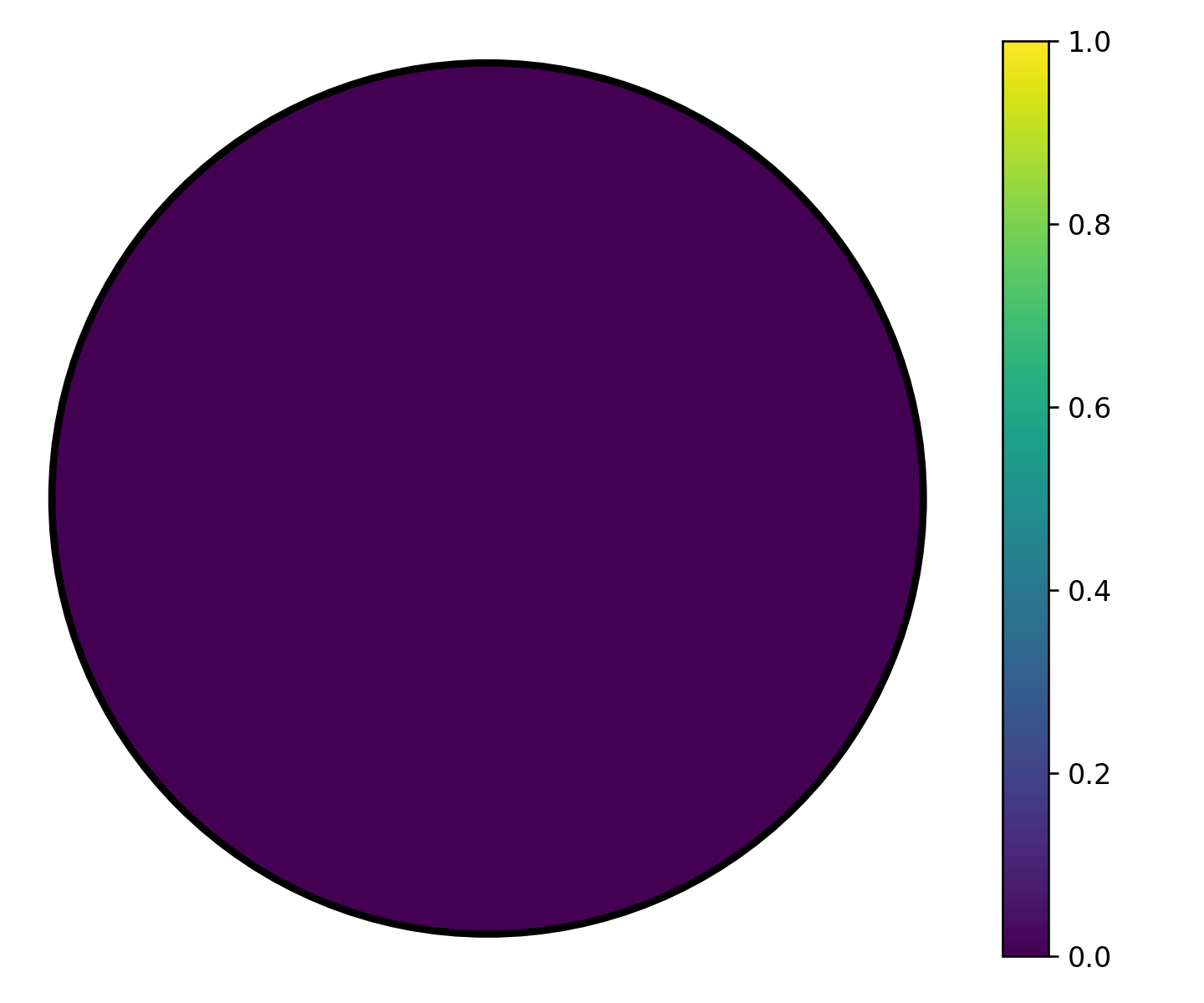}
\end{subfigure}
\caption{Pointwise squared residual $|\tau(u_{\theta})|^{2}$ for a
(perturbed) unknot boundary, before training (left) and after training
(right). Note the colour scales. The residual vanishes along
$\partial D^{2}$ \emph{even before training}: the model is
asymptotically minimal by construction, and the optimiser only works in
the interior. Figure reproduced from \cite{SchettiniUsula2026}.}
\label{fig:before_after}
\end{figure}

Figure \ref{fig:before_after} shows the practical meaning of these
choices: the pointwise residual $|\tau(u_{\theta})|^{2}$ of an
\emph{untrained} model already vanishes along the boundary; training
only needs to remove the error away from the degenerate boundary
region, on compact subsets of the interior where the equation is
uniformly elliptic. In our experience this design is the single largest
contributor to the reliability of the method: with it, plain Adam
\cite{kingma2014adam} followed by L-BFGS \cite{liu1989limited} reaches
expected residuals typically in the $10^{-7}$--$10^{-4}$ range across
the knots we tested (up to $5\times10^{-4}$ for the ten-crossing
$10_{124}$), with essentially no problem-specific tuning.

\begin{rem}
Two comments are needed, for the sake of completeness. 
First, the boundary condition is satisfied exactly
with respect to the degree-$15$ Fourier truncation of $\gamma$ used by
the extension; for the explicit trigonometric parametrisations we use,
the truncation error is negligible, but the statement ``exact boundary
condition'' should be understood in this sense. Second, hard
constraints are not free: they require the constraint manifold
(the space of maps with the prescribed boundary and asymptotics, in this case)
to be explicitly parametrised, which is possible here because the
domain is a disc and the boundary data is a curve. When no such
parametrisation is available, penalty methods remain the fallback.
\end{rem}

\begin{rem}[What did not work]
Since negative results rarely get reported, we list explorations that
were tried and abandoned. Periodic-activation (SIREN-type) correctors
diverged: their large second derivatives, amplified by the $X^{-2}$
factors in \eqref{eq:tau_X}--\eqref{eq:tau_Y}, destabilised training.
Feature embeddings in solid harmonics create a gradient blind spot at
the origin, where all features vanish to high order. Residual-driven
adaptive resampling of collocation points and self-adaptive pointwise
loss weights brought no measurable benefit -- plausibly because the
hard-constrained model leaves no boundary layer for adaptivity to
chase -- and both interact badly with the static-shape compilation
strategy of \S\ref{subsec:compile}. A learnable reparametrisation of
the disc, i.e.~geometric preconditioning, helped in some specific 
cases but was not needed in the final pipeline.
\end{rem}

\section{Evaluating the residual efficiently}\label{sec:efficiency}

The second half of this note concerns a purely computational question:
\emph{what does it cost to evaluate \eqref{eq:loss} and its parameter
gradient, and how should that computation be organised?} The companion
paper records only that second-order derivatives are computed by
composing forward- and reverse-mode automatic differentiation
\cite[\S3.3]{SchettiniUsula2026}, without entering the details about the cost of the
residual evaluation or how to organise it; yet the difference between
the naive and the final organisation is a factor of $40$--$50$ in
wall-clock time (Table \ref{tab:bench}), on the same hardware, for
results identical up to round-off. We describe the two ingredients in
turn.

\subsection{What a training step must compute}\label{subsec:anatomy}

Fix a mini-batch of $N$ collocation points. Evaluating
\eqref{eq:tau_X}--\eqref{eq:LB} at each point requires the value
$u\in\mathbb{R}^{4}$, the Jacobian $J\in\mathbb{R}^{4\times2}$ and the
Hessian $H\in\mathbb{R}^{4\times2\times2}$ of the model with respect to
the \emph{disc coordinates}; everything else -- the pull-back metric
$g=X^{-2}J^{\top}J$, its explicit $2\times2$ inverse, the derivative
terms in \eqref{eq:LB}, the residual itself -- is plain tensor algebra
in $(u,J,H)$: the full residual is an algebraic expression in the
$2$-jet of $u$, with no third derivatives appearing. The optimiser then
needs the gradient of the scalar loss with respect to the
$\sim1.3\times10^{4}$ network parameters.

Automatic differentiation (AD) \cite{baydin2018automatic,griewank2008evaluating}
offers two elementary modes, which a geometer will recognise
immediately: \emph{forward mode} computes pushforwards, mapping an
input tangent vector to an output tangent vector at the cost of one
augmented evaluation; \emph{reverse mode} computes pullbacks of output
covectors. The cost of assembling a full Jacobian is therefore one
sweep per input dimension in forward mode, and one sweep per output
dimension in reverse. For the parameter gradient the choice is forced
and classical: the loss is a single scalar depending on
$\sim10^{4}$ parameters, so reverse mode, i.e.~\textit{backpropagation},
is optimal, and every deep-learning framework provides it.

For the \emph{spatial} derivatives the dimension count points the other
way: the map goes from $\mathbb{R}^{2}$ to $\mathbb{R}^{4}$, so forward
mode needs $2$ sweeps against reverse mode's $4$; and for the Hessian
the disparity compounds. Yet the naive PINN implementation ignores the
count: one calls the framework's reverse-mode routine once per output
component (with the option that keeps the derivative computation itself
differentiable), then differentiates the resulting expressions again
for second derivatives. Our first draft of the code did exactly this, in
the divergence form of \eqref{eq:LB}: per residual evaluation it issued
\emph{sixteen} reverse sweeps -- four for $J$, then twelve more inside
the two divergence-form Laplacians $\Delta_{g}X$ and
$\Delta_{g}\boldsymbol{Y}$ -- each sweep materialising a new
computational graph on top of the previous one, and the whole tower
being differentiated once more, in reverse, for the parameter gradient.
This is the pattern that Table \ref{tab:bench}, row 1, prices at
$427$\,ms per training step.

\subsection{Forward propagation of second-order jets}\label{subsec:jets}

The final implementation computes $(u,J,H)$ in a single forward pass,
with no spatial AD at all, by propagating \emph{second-order jets}
through the network -- forward-mode differentiation carried out
explicitly, order two at a time (in AD terminology, Taylor-mode
\cite{griewank2008evaluating}). An MLP is a composition of affine maps
and coordinatewise activations, and both have trivial jet-transport
rules. Writing $(a,J_{a},H_{a})$ for the value, Jacobian and Hessian of
the signal at the current layer (batched over collocation points):
\begin{align}
\text{affine } a\mapsto Wa+b:\qquad & (a,J_{a},H_{a})\longmapsto(Wa+b,\;WJ_{a},\;WH_{a}),\label{eq:jet_linear}\\
\text{activation } a\mapsto\sigma(a):\qquad & (a,J_{a},H_{a})\longmapsto\bigl(\sigma(a),\;\sigma'(a)\odot J_{a},\;\sigma''(a)\odot J_{a}\otimes J_{a}+\sigma'(a)\odot H_{a}\bigr),\label{eq:jet_act}
\end{align}
with $\odot$ the coordinatewise product and
$(J_{a}\otimes J_{a})_{ij}=\partial_{i}a\,\partial_{j}a$ per neuron;
for $\sigma=\tanh$ one has $\sigma'=1-\sigma^{2}$ and
$\sigma''=-2\sigma\sigma'$, computable from the already-evaluated
activation. Initialising with $J_{a}=\mathrm{Id}_{2}$, $H_{a}=0$ at the
input and applying
\eqref{eq:jet_linear}--\eqref{eq:jet_act} layer by layer yields the
exact $2$-jet of the network at every collocation point in one pass
whose cost is a small constant multiple of a plain forward evaluation
-- roughly the number of jet components carried per neuron, six or
seven depending on whether the symmetry of the Hessian is exploited.

The composite model \eqref{eq:ansatz} multiplies the network by
explicit functions, and jets compose by the Leibniz rule. The $2$-jets
of the two geometric ingredients are available in closed form: for
$\rho_{\mathrm{st}}$ by elementary calculus, and for the
stereobiharmonic extension because differentiating the harmonic
polynomials is algebraic ($\partial_{x}U_{m}=mU_{m-1}$,
$\partial_{x}^{2}U_{m}=m(m-1)U_{m-2}$, and so on), so the jet of the
extension costs one more pass over the same Fourier recursion. The
$2$-jet of \eqref{eq:ansatz} is then assembled by the product rule
-- for instance the first component $X=\rho e^{\nu}$, $\nu=\NN^{X}$,
has
\[
dX=e^{\nu}(d\rho+\rho\,d\nu),\qquad
\mathrm{Hess}\,X=e^{\nu}\bigl(\mathrm{Hess}\,\rho+d\rho\otimes d\nu+d\nu\otimes d\rho+\rho\,d\nu\otimes d\nu+\rho\,\mathrm{Hess}\,\nu\bigr),
\]
and similarly for $Y_{j}=\ext(\gamma)_{j}+\rho^{k}w_{j}$ (with $k$ the
fixed decay exponent of \eqref{eq:ansatz}) -- after which
the residual is evaluated algebraically as in
\S\ref{subsec:anatomy}. Three properties of this organisation matter in
practice.
\begin{enumerate}
\item \emph{It is exact}: all derivatives follow exact rules, and the
result agrees with the nested-autograd implementation to
$5\times10^{-14}$ in double precision (\S\ref{subsec:measured}), i.e.\
to accumulated round-off.
\item \emph{It separates the two differentiations}: spatial derivatives
no longer pass through the AD engine at all -- collocation points are
plain data, and no graph of graphs is ever built; the single remaining
reverse pass, for the parameter gradient, differentiates an ordinary
composition of tensor operations.
\item \emph{It mirrors the mathematics}: the implementation is
literally the chain rule, the Leibniz rule and the coordinate formula
\eqref{eq:LB}, in the order in which one would write them on paper.
This makes it auditable line by line -- a non-trivial virtue when the
output feeds a conjecture test -- and it makes the whole residual a
\emph{flat, branch-free tensor program}, which is precisely what the
compiler discussed next requires.
\end{enumerate}
On its own, this reorganisation accounts for a factor $\approx22$--$24$
(Table \ref{tab:bench}, row 3).

A remark on precision, complementing \cite[\S3.3]{SchettiniUsula2026}:
we run everything in double precision, since in single precision the
rounding incurred in second-derivative computations makes the residual
saturate near $10^{-4}$ even at a true solution, defeating both the
L-BFGS refinement and the double-point analysis. This choice interacts
with hardware: consumer GPUs execute \texttt{float64} at a small
fraction of their nominal throughput (and Apple-silicon GPUs not at
all), so the pipeline was designed \emph{CPU-first} -- which raises the
stakes for the compilation strategy below, as the CPU has no raw-power
headroom to hide overheads in.

\subsection{Compiling the computational graph once}\label{subsec:compile}

A PyTorch program normally runs in \emph{eager mode}: every tensor
operation is dispatched to a pre-compiled kernel as the Python
interpreter encounters it, and the computational graph exists only
implicitly, rebuilt at every step for the benefit of the backward pass.
For a residual like ours, i.e.~hundreds of small operations on modest
batches, the per-operation overhead (Python dispatch, kernel launch,
allocation of intermediates) dominates the arithmetic. The remedy,
available in modern frameworks under names such as
\texttt{torch.compile} \cite{ansel2024pytorch2} or \texttt{jax.jit},
is to \emph{trace} the computation once into an explicit graph,
optimise it (fusing chains of elementwise operations into single
kernels, eliminating dead code and redundant intermediates), generate
machine code -- C++ kernels, in our CPU case -- and thereafter re-run
the compiled artefact, forward and backward, at every step.

The catch is that tracing specialises: the compiled artefact is valid
for a fixed operation sequence and (in the static regime we use) fixed
tensor shapes, and anything that varies from step to step triggers
either a fallback to eager execution or an expensive silent
recompilation. Concretely, the final implementation:
\begin{enumerate}
\item packages the entire residual evaluation --- jet propagation
\eqref{eq:jet_linear}--\eqref{eq:jet_act}, closed-form jets of $\rho$
and $\ext(\gamma)$, Leibniz assembly, metric algebra, residual --- as
one flat module, compiled \emph{once} as a single graph with static
shapes (\texttt{fullgraph=True}, \texttt{dynamic=False}); the Fourier
coefficients of the extension enter as constant buffers baked into the
graph;
\item makes the training loop shape-static to match: the collocation
pool ($2^{14}$ points) is drawn once and only re-shuffled between
epochs, mini-batches have a fixed size ($2^{10}$; a remainder that
would produce a ragged final batch is dropped), and the points carry no
AD metadata;
\item compiles one further instance for the L-BFGS refinement phase,
whose full-batch closure has a different (but again fixed) input shape;
each compiled graph is reused for the entire phase -- up to $1.6\times10^{5}$
Adam steps and $10^{4}$ full-batch L-BFGS iterations in a production
run.
\end{enumerate}
The one-off cost of compilation is $14$--$18$ seconds in our setting; 
amortised across a training run it is negligible, but it is
\emph{per shape}: an innocently varying batch size would pay it again
and again, which is the practical reason for the static sampling
design in item (2). It also explains a symptom worth knowing:
a shape change mid-training manifests as a multi-second,
log-silent stall that is easily mistaken for a hang.

We add one experience report. Before settling on hand-rolled jets we
used the framework's composable functional transforms
(\texttt{torch.func}'s \texttt{jacfwd}/\texttt{jacrev}) to obtain
$(u,J,H)$; these already realise the forward-mode dimension count and
run an order of magnitude faster than nested autograd (Table
\ref{tab:bench}, row 2); for some time, they were the production
pipeline. But composing them \emph{with the graph compiler} proved
fragile: in our environment (PyTorch 2.9), higher-order functional
transforms over a function containing non-trivial input arithmetic
could not be traced. The workaround was structural: replace
transform-generated derivatives by the explicit jets of
\S\ref{subsec:jets}, so that the program presented to the compiler
contains only elementary tensor operations. Relative to this eager
transform pipeline (the state of the code immediately before
compilation was introduced), the compiled jet evaluator is
$2.9$--$4.9\times$ faster (Table \ref{tab:bench}, rows 2 and 4): this,
rather than the smaller margin over the eager jets, is the improvement
that was actually experienced when the compiled graph entered the
pipeline. The general lesson is that the \emph{intersection} of
advanced features (higher-order AD, functional transforms, compilers)
is the least-charted part of any framework, and computations organised
as plain tensor programs age better than ones that lean on that
intersection.

\subsection{Measured impact}\label{subsec:measured}

Table \ref{tab:bench} isolates the two contributions on the production
task of \cite{SchettiniUsula2026}: identical model (perturbed-trefoil
boundary, default architecture, fixed random initialisation -- the
timings are insensitive to the weights), identical collocation points,
identical hardware, double precision throughout; each timed step
evaluates the residual on the batch, forms the mean-square loss, and
computes the full parameter gradient, i.e.~the work of one optimisation
step, apart from the parameter update itself, for the two batch sizes
used by the two training phases. The baseline row runs the project's
original nested reverse-mode implementation; before timing, we
verified that all four pipelines return the same residual field, with
maximum pairwise deviation $5.0\times10^{-14}$ on a residual field of
maximum modulus $\approx53$.

\begin{table}[t]
\centering
\begin{tabular}{lrrrr}
\toprule
 & \multicolumn{2}{c}{$N=2^{10}$ (Adam)} & \multicolumn{2}{c}{$N=2^{14}$ (L-BFGS)}\tabularnewline
Residual pipeline & ms/step & speedup & ms/step & speedup\tabularnewline
\midrule
nested reverse-mode AD (\S\ref{subsec:anatomy}) & $427.5$ & $1$ & $3260.7$ & $1$\tabularnewline
\texttt{torch.func} transforms, eager (\S\ref{subsec:compile}) & $42.0$ & $10.2\times$ & $223.7$ & $14.6\times$\tabularnewline
forward $2$-jets, eager (\S\ref{subsec:jets}) & $17.6$ & $24.3\times$ & $143.9$ & $22.7\times$\tabularnewline
forward $2$-jets, compiled (\S\ref{subsec:compile}) & $8.6$ & $49.7\times$ & $78.2$ & $41.7\times$\tabularnewline
\bottomrule
\end{tabular}
\medskip
\caption{Median wall-clock time per training step (residual + loss +
parameter gradient) over $20$ steps after warm-up, on a laptop CPU
(Apple M4 Pro, $10$ threads, PyTorch 2.9.1, Python 3.12,
\texttt{float64}); speedups are computed from the unrounded timings.
The \texttt{torch.func} row reports the faster of the
forward-over-reverse and forward-over-forward compositions at each
batch size. One-off graph compilation costs $14$--$18$\,s per input
shape (excluded: a one-off cost at the start of each training phase).
The four pipelines agree pairwise to $5\times10^{-14}$ on the same
inputs.}
\label{tab:bench}
\end{table}

Let us conclude with three observations. First, the factors compose as
claimed: moving from nested reverse-mode AD to the forward-mode
functional transforms is worth $10$--$15\times$, replacing the
transforms by explicit jets a further $1.6$--$2.4\times$, and
compiling the resulting graph a further $1.8$--$2.1\times$, for
$42$--$50\times$ end to end. Second, the gains are broadly stable across a $16$-fold change in batch size,
so they reflect the algorithm, not a small-batch overhead artefact.
Third, the absolute numbers set the
scale of experiments one can contemplate; at the production schedule
($10^{4}$ Adam epochs of $16$ mini-batches plus up to $10^{4}$
full-batch L-BFGS iterations), the compiled pipeline completes the
optimisation in roughly $40$ minutes of compute on a laptop (consistent
with the roughly one hour reported in the companion paper) whereas the same
schedule through the baseline pipeline extrapolates to more than a
day. The entire study of \cite{SchettiniUsula2026}, some dozens of
trained models plus failed runs and ablations, was carried out on
ordinary laptops; without these two optimisations it would have
required a cluster, or would not have been done.

\section{Guidelines}\label{sec:guidelines}

We summarise the discussion as advice to a geometric analyst or differential geometer starting a
PINN project. None of it is specific to minimal surfaces.

\begin{enumerate}
\item \emph{Parametrise constraints away; penalise only what you
cannot parametrise.} Build boundary conditions, positivity, decay and
symmetry into the model: every constraint moved from the loss into the
architecture removes a weighting hyperparameter, shrinks the search
space, and holds identically throughout training.
\item \emph{Spend your regularity theory before you spend your
compute.} Whatever is known about solutions a priori (boundary
asymptotics, decay rates, symmetries) is structure the network would
otherwise have to learn, imperfectly, from the residual signal alone.
\item \emph{Make $\theta=0$ meaningful.} Arrange the architecture so
the zero network is the best explicit approximation you can write down,
and check that this initial map is smooth -- the network will not fix
a singular ansatz.
\item \emph{Count dimensions before differentiating.} Reverse mode for
the parameter gradient; forward mode for spatial derivatives whenever
the domain dimension is small, as it is for curves and surfaces. For
second-order residuals, propagating the $2$-jet explicitly through the
network --- a few dozen lines of code implementing equations
\eqref{eq:jet_linear}--\eqref{eq:jet_act} -- is exact, fast, and
auditable.
\item \emph{Compile once; design for it.} Keep shapes static (fixed
collocation pool, fixed batch size), keep the residual a flat
branch-free tensor program, and expect one compiled artefact per input
shape. Treat a mysterious mid-run stall as a recompilation until proven
otherwise.
\item \emph{Use double precision, and measure.} Second derivatives in
single precision floor the residual near $10^{-4}$; \texttt{float64} on
a compiled CPU pipeline is a perfectly viable regime. Verify optimised
pipelines against a naive reference implementation to round-off, and
time the actual training step, not a proxy.
\end{enumerate}

\section{Concluding remarks}\label{sec:conclusion}

The asymptotic Plateau problem is, in retrospect, an ideal proving
ground for geometry-aware PINNs: the boundary condition sits at
infinity, the boundary asymptotics of solutions are theorems, an exact
model solution exists to anchor the ansatz, and the probed
conjecture (signed counts of double points confronted with HOMFLY
coefficients \cite{FineKnots}) demands a precision that soft
constraints and single precision cannot deliver. But nothing in the two
theses of this note is specific to it: hard-constrained architectures
apply whenever the constraint set can be parametrised, and the
jet-propagation and compilation techniques of \S\ref{sec:efficiency}
apply verbatim to any PINN whose domain has low dimension -- the
typical situation for parametrised curves, surfaces and maps in
geometric analysis. The efficiency margin of Table \ref{tab:bench} is
also what makes the natural next steps thinkable at
laptop-to-workstation scale: higher-genus domains, systematic knot
families, and, most importantly, pushing residuals towards the regime where the numerical solutions could serve as the
starting point of computer-assisted existence proofs, in the spirit of
\cite{gomez2019computer,wang2025discoveryunstablesingularities,platt2026nonuniquenesssymmetriesnirenbergproblem}.
We hope this note lowers the activation energy for colleagues in
geometric analysis and differential geometry to attempt such computations themselves.

\begin{small}
\bibliographystyle{plain}
\bibliography{Bibliography}
\end{small}

\end{document}